# Second Order Step by Step Sliding mode Observer for Fault Estimation in a Class of Nonlinear Fractional Order Systems

Seyed Mohammad Moein Mousavi
Student
Electrical and Computer Engineering Department
Tarbiat Modares University
Tehran, Iran
Email: moein_mousavi@modares.ac.ir

Amin Ramezani
Professor
Electrical and Computer Engineering Department
Tarbiat Modares University
Tehran, Iran
Email: ramezani@modares.ac.ir

*Abstract*—**This paper considers fault estimation in nonlinear fractional order systems in observer form. For this aim, a step by step second order sliding mode observer is used. By means of a fractional inequality, the stability of the observer estimation errors is analyzed and some conditions are introduced to guarantee finite time convergence of estimation errors. Finally, in a numerical example, effectiveness of this observer is demonstrated.**

## 1 INTRODUCTION

Fractional calculus as a generalization of classic calculus is a mathematic tool which is recently being used in control engineering. Among recent years it is proven that some systems are modeled more accurate using fractional order models, compare to integer order ones [2].

Also fault diagnosis has been always an important topic in the industry [4]. Faults can cause failure and damage in physical systems if they are not detected in appropriate time. Fault detection methods are divided into two main types: data based and model based methods. Indeed, model based fault detection methods have been extended during last three decades. In these methods, sensor and actuator faults are detected through the relations between accessible signals. the most popular model based methods are: parameter estimation, observer design [5] and parity space. After a fault is detected, it is important to estimate its shape and domain as an unknown input. Since some systems are modeled by fractional models, using fractional order observers for fault estimation is an important issue.

In [6], observability of the states in nonlinear fractional order systems is discussed and using a first order sliding mode observer, fault is estimated as an unknown input in such system. But as the observer is a first order one, the chattering problem exists. Sliding mode observer for fault estimation is discussed in [13] and unknown input observer for fractional order system is considered in [10].in [12] state estimation in a nonlinear fractional order system with uncertain parameters using sliding mode observers is taken into account. sliding mode controller design for a nonlinear fractional order system for disturbance rejection is discussed in [9]. Also in [7,8,11] the structure of a second order sliding mode observer for a nonlinear integer order system is introduced for fault estimation.in [14] a second order sliding mode observer is used for fault estimation in a linear fractional order system.

Hence to the best of our knowledge, fault estimation in nonlinear fractional order systems using second order sliding mode observer (which is less affected by the chattering problem) is not discussed in the literature. However, an unknown input observer for this type of systems is designed in [15] which does not consider the unknown input as a state in the observer structure and it is going to be discussed is this paper. In the following parts, some preliminaries on fractional calculus is introduced in section 2. The observer structure and conditions for stability of the state estimation error is introduced in section 3.in section 4 the effectiveness of the proposed observer is discussed by a numerical example and simulation. and finally in section 5 we have the conclusion part.

## 2 PRELIMINARIES

Let $C[a\ b]$ be the space of continuous functions $f(t)$ on $[a\ b]$ and we mean by $C^k$ the space of real-valued functions $f(t)$ with continuous derivatives up to order $k-1$ such that $f^{(k-1)}(t) \in C[a\ b]$ and $f^i(t)$ is the i-th derivative of $f(t)$.

### 2.1 fractional calculus

According to [1] there are three main definitions of fractional order derivatives:

- The Riemann–Liouville fractional derivative of order $\alpha$ of $f(t) \in C^m[a\ b]; t \in [a\ b]$:

$$^{RL}_aD^\alpha_t f(t) = \frac{1}{\Gamma(m-\alpha)} \frac{d^m}{dt^m} \int_a^t (t-\tau)^{m-\alpha-1} f(\tau) d\tau \qquad (1)$$

- Caputo's derivative of order $\alpha$ of $f(t) \in C^m[a\ b]; t \in [a\ b]$:

$${}^{C}_{a}D^{\alpha}_{t}f(t) = \frac{1}{\Gamma(m-\alpha)}\int_{a}^{t}(t-\tau)^{m-\alpha-1}\frac{d^m f(\tau)}{d\tau^m}d\tau \quad (2)$$

- Grunwald-Letnikov definition:

$${}^{GL}_{a}D^{\alpha}_{t}f(t) = \lim_{h\to\infty}\frac{1}{\Gamma(\alpha)h^{\alpha}}\sum_{j=0}^{\left[\frac{t-a}{h}\right]}\frac{\Gamma(\alpha+j)}{\Gamma(j+1)}f(t-jh) \quad (3)$$

Where $\left[\frac{t-a}{h}\right]$ denotes the integer part of $\frac{t-a}{h}$ and $\Gamma()$ is Euler's gamma function.

The drawback of the first definition is that the initial conditions are in terms of the variable's fractional order derivatives. having said that, Caputo's definition of fractional derivative needs the initial conditions of the main function and not its fractional derivatives. hence, in engineering usages caputo definition is commonly applied. for the simplicity in the rest of this paper this notation is used:

$${}^{C}_{0}D^{\alpha}_{t}f(t) \triangleq D^{\alpha}f(t) \quad (4)$$

*2.2 Fractional inequality*

This lemma is proven in [3] which will help us with stability analyze in fractional order systems:

**Lemma1**: Let $x(t) \in^n$ be a differentiable vector. For any $t \geq t_0$:

$$\frac{1}{2}D^{\alpha}(e^T P e) \leq e^T P D^{\alpha}e \quad (5)$$

Where $P \in^{n \times n}$ is a symmetric positive definite and constant matrix

## 3 OBSERVER DESIGN

*3.1 Observer structure*

Let us consider the system in this observable form:

$$\begin{cases} D^{\alpha}x_1 = x_2, \ldots, D^{\alpha}x_{n-1} = x_n \\ D^{\alpha}x_n = f_1(x) + f_2(x)f \end{cases} \quad (6)$$

Where $0 < \alpha < 1$ and $x \in^{n \times 1}$ is the state vector, $f_1(x)$ and $f_2(x)$ are Lipchitz functions and all the states are assumed to be bounded:

$$\text{for all } t > 0 : |x_i(t)| < a_i \quad (7)$$

And for all $t > 0$:

$$|f(t)| < A_1$$
$$|f_1(t)| < A_2$$
$$|f_2(t)| < A_3$$
$$|D^{\alpha}f(t)| < \acute{A}_1$$
$$|D^{\alpha}f_1(t)| < \acute{A}_2$$
$$|D^{\alpha}f_2(t)| < \acute{A}_3$$

$A_i$ and $\acute{A}_i$ are known positive constants. the observer structure is considered as:

$$\begin{cases} D^{\alpha}\hat{x}_1 = \tilde{x}_2 + \lambda_1|e_1|^{0.5}sign(e_1) \\ D^{\alpha}\tilde{x}_2 = \alpha_1 sign(e_1) \\ D^{\alpha}\hat{x}_2 = E_1[\tilde{x}_3 + \lambda_2|e_2|^{0.5}sign(e_2)] \\ D^{\alpha}\tilde{x}_3 = E_1\alpha_2 sign(e_2) \\ D^{\alpha}\hat{x}_3 = E_2[\tilde{x}_4 + \lambda_3|e_3|^{0.5}sign(e_3)] \\ \vdots \\ D^{\alpha}\tilde{x}_{n-1} = E_{n-3}\alpha_{n-2}sign(e_{n-2}) \\ D^{\alpha}\hat{x}_{n-1} = E_{n-2}[\tilde{x}_n + \lambda_{n-1}|e_{n-1}|^{0.5}sign(e_{n-1})] \\ D^{\alpha}\tilde{x}_n = E_{n-2}\alpha_{n-1}sign(e_{n-1}) \\ D^{\alpha}\hat{x}_n = E_{n-1}[f_1(\tilde{x}) + f_2(\tilde{x})\tilde{f} + \lambda_n|e_n|^{0.5}sign(e_n)] \\ D^{\alpha}\tilde{f} = E_{n-1}\alpha_n sign(e_n) \\ D^{\alpha}\hat{f} = E_n\left[\tilde{\theta} + \lambda_{n+1}|e_f|^{0.5}sign(e_f)\right] \\ D^{\alpha}\tilde{\theta} = E_n\alpha_{n+1}sign(e_f) \end{cases} \quad (8)$$

Considering fault as a state in the observer structure causes more complexity and the more number of gains to be chosen. This will make the design harder but the estimations will be more accurate. Estimation errors are defined as:

$$e_i = \tilde{x}_i - \hat{x}_i, e_f = \tilde{f} - \hat{f} \quad (9)$$

With $\tilde{x}_1 = x_1$ and $x_1$ is the output and the only accessible signal. We have:

$$E_i = 1 \text{ if } |e_j| = |\tilde{x}_j - \hat{x}_j| \leq \varepsilon \text{ for all } j \leq i \text{ else } E_i = 0$$

$\varepsilon$ is a small positive scalar and the observer gains are all positive coefficients.

*3.2 Stability analyze*

Step1: $E_i = 0 \text{ for all } i$

using (6) and (8), error dynamics are obtained as:

$$\begin{cases} D^{\alpha}e_1 = x_2 - \tilde{x}_2 - \lambda_1|e_1|^{0.5}sign(e_1) \\ D^{\alpha}\tilde{x}_2 = \alpha_1 sign(e_1) \end{cases} \quad (10)$$

Using lemma 1:

$$D^{2\alpha}e_1 \leq x_3 - \alpha_1 sign(e_1)$$
$$-\frac{1}{2}\lambda_1 D^{\alpha}e_1|e_1|^{-0.5} \quad (11)$$

If the following equations hold, $e_1$ dynamic would be stable and it will converge to zero in bounded time.

$$\alpha_1 > \alpha_3$$
$$\lambda_1^2 > \frac{4\alpha_3(\alpha_1 + \alpha_3)}{\alpha_1 - \alpha_3}$$

So $e_1$ and its fractional derivatives would become zero. Considering (10), $x_2$ will equal to $\tilde{x}_2$.

Step 2: $E_1 = 1$ and:
$$D^\alpha e_1 = 0 \quad (12)$$
$$D^\alpha e_2 = x_3 - \tilde{x}_3 - \lambda_2 |e_2|^{0.5} sign(e_2) \quad (13)$$
$$D^\alpha \tilde{x}_3 = \alpha_2 sign(e_2) \quad (14)$$

Similarly, with a suitable choice of observer gains, $e_2$ and its derivatives would be 0 and $x_3$ will equal to $\tilde{x}_3$.

Step n: $E_1 = E_2 = \cdots = E_{n-1} = 1$ and
$$D^\alpha e_1 = \cdots = D^\alpha e_{n-1} = 0 \quad (15)$$
$$D^\alpha e_n = f_1(x) + f_2(x)f \\ -f_1(\tilde{x}) - f_2(\tilde{x})\tilde{f} - \lambda_n |e_n|^{0.5} sign(e_n) \quad (16)$$
$$D^\alpha \tilde{f} = \alpha_n sign(e_n) \quad (17)$$

Again, using lemma1 we have:
$$D^{2\alpha} e_n \leq D^\alpha (f_1(x) - f_1(\tilde{x})) + \\ (D^\alpha f_2(x))f + f_2(x)(D^\alpha f) - (D^\alpha f_2(\tilde{x}))\tilde{f} \\ -f_2(\tilde{x})(D^\alpha \tilde{f}) - \frac{1}{2}\lambda_n D^\alpha e_n |e_n|^{-0.5} \quad (18)$$

Because of a similar reason as it is assumed that all the states, fault and their derivatives are bounded, there exists some values for observer gains such that dynamic of $e_n$ becomes stable so $\tilde{f}$ converges to $f$ and again, if the observer gains are chosen big enough, all the estimations converge to states. Hence the fault signal will be reconstructed.

## 4 NUMERICAL SIMULATION

Consider a nonlinear fractional order system as:
$$\begin{cases} D^\alpha x_1 = x_2 \\ D^\alpha x_2 = x_3 \\ D^\alpha x_3 = -0.5 x_1 - \sin(x_2) - x_3 |x_3| + f \\ f = 0.5 \cos(0.5\pi t) \end{cases}$$

$f(t)$ denotes the fault to be estimated, $\alpha = 0.7$ And the observer structure is:
$$\begin{cases} D^\alpha \hat{x}_1 = \tilde{x}_2 + \lambda_1 |e_1|^{0.5} sign(e_1) \\ D^\alpha \tilde{x}_2 = \alpha_1 sign(e_1) \\ D^\alpha \hat{x}_2 = E_1[\tilde{x}_3 + \lambda_2 |e_2|^{0.5} sign(e_2)] \\ D^\alpha \tilde{x}_3 = E_1 \alpha_2 sign(e_2) \\ D^\alpha \hat{x}_3 = E_2[f_1(\tilde{x}) + f_2(\tilde{x})\tilde{f} + \lambda_3 |e_3|^{0.5} sign(e_3)] \\ D^\alpha \tilde{f} = E_2 \alpha_3 sign(e_3) \\ D^\alpha \hat{f} = E_3 \left[\tilde{\theta} + \lambda_4 |e_f|^{0.5} sign(e_f)\right] \\ D^\alpha \tilde{\theta} = E_3 \alpha_4 sign(e_f) \end{cases}$$

$x_1$ as the output of the system is the only accessible signal and it is the observer input. The initial conditions of the observer parameters are zero and $x_{1_0} = 0.1$, $x_{2_0} = 0.1$, $x_{3_0} = -0.1$. Observer gains are $\lambda_1 = \lambda_2 = \lambda_3 = \lambda_4 = 0.1, \alpha_1 = 1, \alpha_2 = 2, \alpha_3 = 5, \alpha_4 = 10$. It should be noted that observer gains are chosen by trial and error. If these gains are chosen small then the estimation errors will not converge to zero, also if we choose big values for the gains, the chattering problem will appear. Therefore, a future work can be optimum choice of these gains.

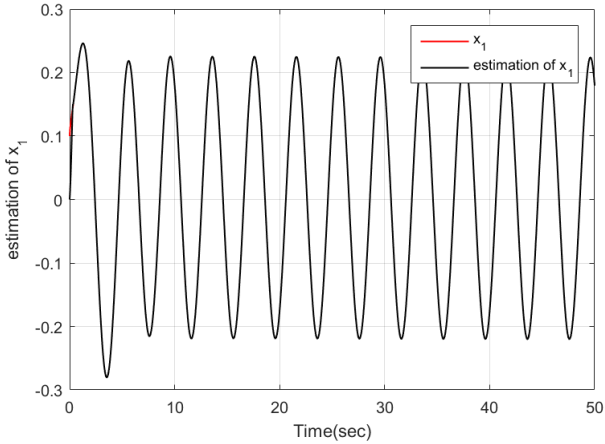

Figure 1. estimation of $x_1$

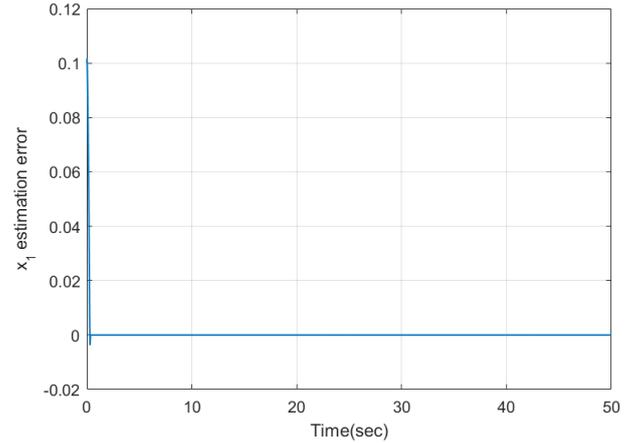

Figure 2. estimation error of $x_1$

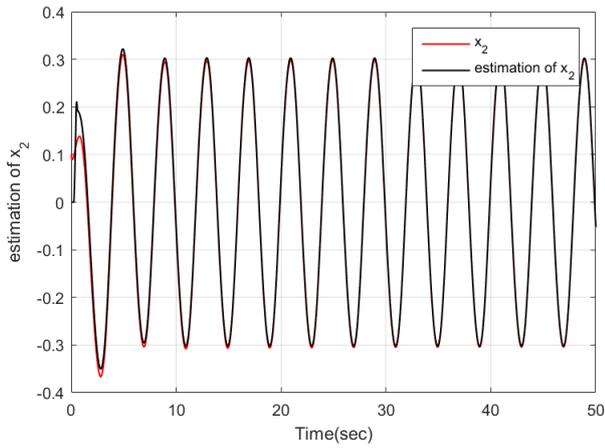

Figure 3. estimation of $x_2$

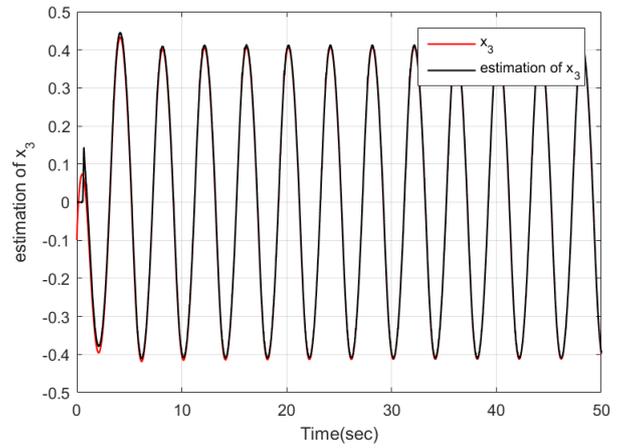

Figure 5. estimation of $x_3$

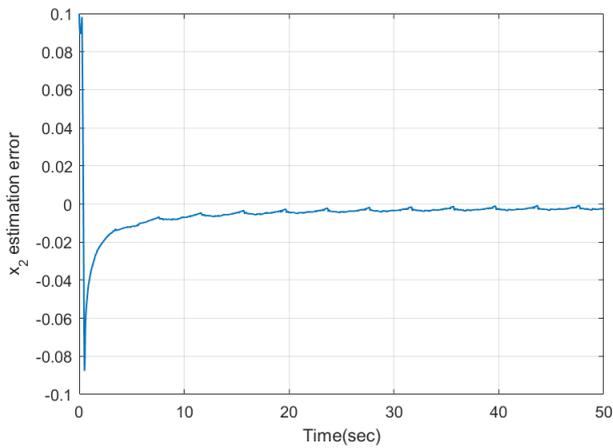

Figure 4. estimation error of $x_2$

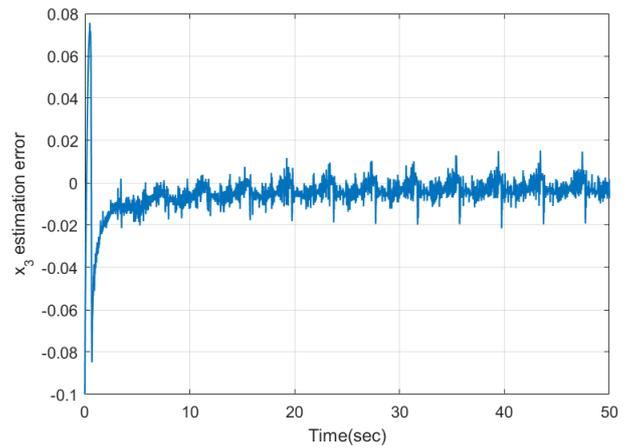

Figure 6. estimation error of $x_3$

Simulation Is performed using Grünwald–Letkinov's definition of the fractional derivative. State estimations and their errors are depicted in Figures.1-6. Also Figures.7-8 represent the fault signal and its estimation and the estimation error. It is clear that the fault signal is reconstructed and the estimation errors converge to zero. As a second order sliding mode observer is adopted here, it is less affected by the chattering problem compare to first order case. Hence this simulation endorses the effectiveness of the proposed observer.

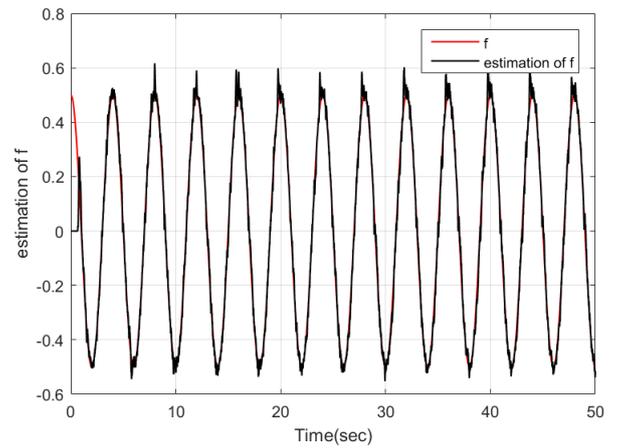

Figure 7. estimation of $f$

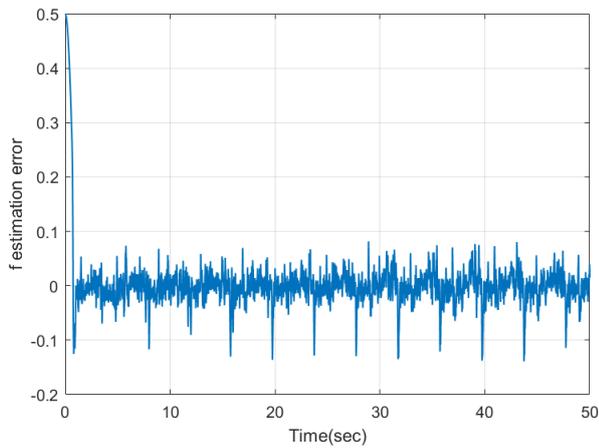

Figure 8. estimation error of $f$

## 5 CONCLUSION

In this paper a second order step by step sliding mode observer was adopted to estimate states and fault as an unknown input in a nonlinear fractional order system. Using a fractional inequality, conditions for stability and finite time convergence of estimation errors were introduced. a numerical example of a commensurate 3-dimensional nonlinear fractional order system was represented to illustrate the effectiveness of the proposed observer. Since there is no specific method to select observer gains, a future work on this topic can be optimum choice of gains in order to minimize the estimation error.